\documentclass{amsart}
\usepackage{latexsym}
\usepackage{amssymb,amsmath,amsfonts}

\newtheorem{theorem}{Theorem}

\newtheorem{lemma}[theorem]{Lemma}

\theoremstyle{definition}
\newtheorem{definition}[theorem]{Definition}

\newtheorem{example}[theorem]{Example}
\newtheorem{rmk}[theorem]{Remark}
\theoremstyle{remark}

\begin{document}

\title[Corrigendum] {\bf Corrigendum to the paper, ``A new iteration process for approximation of common fixed points for finite families of total asymtotically nonexpansive mappings". \\{\sf int. J. Math. Math. Sci.} vol. 2009, doi:10.1155/2009/615107.}

\author{C. E. Chidume$^{1}$}
\address{{$~^{1}$}The Abdus Salam International Centre for Theoretical
Physics, Trieste, Italy}\email{chidume@ictp.trieste.it}
\author{E.U. Ofoedu$^{2}$\\
  }
\address{{$~^{2}$} Department of Mathematics, Nnamdi Azikiwe University, P.M.B. 5025, Awka,
Anambra State, Nigeria} \email{euofoedu@yahoo.com}

\keywords{{uniformly $L$- Lipschitzian maps, asymptotically
nonexpansive mappings, total asymptotically quasi-nonexpansive
mappings, Banach spaces.
}\\
 {$^{2}$}This author undertook this work when he was atthe Abdus Salam International Centre for Theoretical Physics
 as a visiting fellow.}

\begin{abstract} A gap in the proof of Theorem 3.5 in the above paper is observed. The argument used on page 11, starting from line 8 from bottom to the end of the proof of the theorem is not correct. In this corrigendum, it is our aim to close this gap.
\end{abstract}

\maketitle \thispagestyle{empty}

\section{Introduction}

\noindent Let $K$ be a nonempty subset of a real normed 
space $E.$ A mapping $T:K\to K$ is said to be \emph{nonexpansive}
if $\|Tx-Ty\|\le\|x-y\|$ for all $x,y\in K.$

\noindent The mapping $T$ is called \emph{asymptotically
nonexpansive} if there exists a sequence $\{\mu_n\}_{n\ge1}\subset
[0, \infty)$ with $\displaystyle \lim_{n\to \infty}\mu_n=0$ such
that for all $x,y\in K,$
\begin{eqnarray*}
\|T^nx-T^ny\|\le(1+\mu_n)\|x-y\|\;for \;all\;n\ge 1;
\end{eqnarray*} and $T$ is said to be uniformly $L$-Lipschitzian if there exists
a constant $L\ge 0$ such that 
$$\|T^n x- T^n y\|\le L\|x-y\|\;\forall\;x,y\in K.$$

\noindent The class of asymptotically nonexpansive mappings was
introduced by Goebel and Kirk \cite{Goebel} as a generalization of
the class of nonexpansive mappings. They proved that if $K$ is a
nonempty closed convex bounded subset of a uniformly convex real 
Banach space and $T$ is an asymptotically nonexpansive
self-mapping of
$K,$ then $T$ has a fixed point.\\

\noindent A mapping $T$ is said to be \emph{asymptotically
nonexpansive in the intermediate sense} (see e.g., \cite{Bruck})
if it is continuous and the following inequality holds:
\begin{eqnarray}\label{e1}
\displaystyle\limsup_{n\to\infty}\sup_{x,y\in
K}(\|T^nx-T^ny\|-\|x-y\|)\le 0.
\end{eqnarray}
Observe that if we define
\begin{eqnarray*}
\displaystyle a_n:=\sup_{x,y\in K}(\|T^nx-T^ny\|-\|x-y\|),\;and\;
\sigma_n=max\{0, a_n\},
\end{eqnarray*}
then $\sigma_n\to 0$ as $n\to\infty$ and \eqref{e1} reduces to
\begin{eqnarray}\label{e3}
\|T^nx-T^ny\|\le\|x-y\|+\sigma_n,\;for \;all\;x,y\in K,\;n\ge 1.
\end{eqnarray}

\noindent The class of mappings which are asymptotically
nonexpansive in the intermediate sense was introduced by Bruck
\emph{et al.} \cite{Bruck}. It is known \cite{Kirk} that if $K$ is
a nonempty closed convex bounded subset of a uniformly convex real 
Banach space $E$ and $T$ is a self-mapping of $K$ which is
asymptotically nonexpansive in the intermediate sense, then $T$
has a fixed point. It is worth mentioning that the class of
mappings which are asymptotically nonexpansive in the intermediate
sense contains properly the
class of asymptotically nonexpansive mappings (see, e.g., \cite{Kim}).\\

\noindent Sahu \cite{Sahu}, introduced the class of nearly Lipschitzian 
mappings. Let $K$ be a nonempty subset of a normed space $E$ and let $\{a_n\}_{n\ge 1}$
be a sequence in $[0,+\infty)$ such that $\displaystyle\lim_{n\to\infty}a_n=0.$ A mapping 
$T:K\to K$ is called \emph{nearly Lipschitzian} with respect to $\{a_n\}_{n\ge 1}$ if for 
each $n\in \mathbb{N},$ there exists $k_n\ge 0$ such that 
\begin{eqnarray}\label{en}
\|T^nx-T^ny\|\le k_n(\|x-y\|+a_n)\;\forall\;x,y\in K.
\end{eqnarray} Define $$\eta(T^n):=\sup\Bigl\{\frac{\|T^nx-T^ny\|}{\|x-y\|+a_n}:x,y\in K,\;x\ne y\Bigr\}.$$
Observe that for any sequence $\{k_n\}_{n\ge 1}$ satisfying \eqref{en}, $\eta(T^n)\le k_n\;\forall\;n\in \mathbb{N}$
and that $$\|T^nx-T^ny\|\le \eta(T^n)(\|x-y\|+a_n)\;\forall\;x,y\in K,\;n\in \mathbb{N}.$$ $\eta(T^n)$ is called the 
\emph{nearly Lipschitz constant} of the mapping $T$. A nearly Lipschitzian mapping $T$ is said to be
\begin{itemize}
\item \emph{nearly contraction} if $\eta(T^n) < 1$ for all $n\in\mathbb{N};$
\item \emph{nearly nonexpansive} if $\eta(T^n)=1$ for all $n\in \mathbb{N};$
\item \emph{nearly asymptotically nonexpansive} if $\eta(T^n)\ge 1$ for all $n\in\mathbb{N}$ and \\
$\displaystyle\lim_{n\to\infty}\eta(T^n)=1;$
\item \emph{nearly uniform $L$-Lipschitzian} if $\eta(T^n)\le L$ for all $n\in \mathbb{N};$
\item \emph{nearly uniform $k$-contraction} if $\eta(T^n)\le k < 1$ for all $n\in \mathbb{N}.$
\end{itemize}

\begin{example}(See Sahu \cite{Sahu})
\noindent Let $E=\mathbb{R},\;K=[0,1].$ Define $T:K\to K$ by 
\begin{eqnarray*}
Tx={} \left\{ \begin{array} {ll} \frac{1}{2},\;if\;x\in [0,\frac{1}{2}]\\
0,\;if \;x\in (\frac{1}{2},1].\\
\end{array}\right.
\end{eqnarray*} It is obvious that $T$ is not continuous, and thus,
 not Lipschitz. However, $T$ is nearly nonexpansive. In fact, for a real
 sequence $\{a_n\}_{n\ge 1}$ with $a_1=\frac{1}{2}$ and $a_n\to 0$ as $n\to\infty,$
 we have $$\|Tx-Ty\|\le \|x-y\|+a_1\;\forall\;x,y\in K$$ and 
 $$\|T^nx-T^ny\|\le \|x-y\|+a_n\;\forall\;x,y\in K,\;n\ge 2.$$ 
 This is becuase $T^nx=\frac{1}{2}\;\forall\;x\in [0,1],\;n\ge 2.$
\end{example}

\begin{rmk}\label{rmk1}
If $K$ is a bounded domain of an asymptotically nonexpansive mapping $T,$ 
then $T$ is nearly nonexpansive. In fact, for all $x,y\in K$ and $n\in\mathbb{N},$
 we have $$\|T^nx-T^ny\|\le (1+\mu_n)\|x-y\|\le \|x-y\|+ diam(K)\mu_n.$$
 Furthermore, we easily observe that every nearly nonexpansive mapping 
 is nearly asymptotically nonexpansive with $\eta(T^n)\equiv 1\;\forall\;n\in \mathbb{N}.$
\end{rmk}

\begin{rmk}\label{rmk2}
If $K$ is a bounded domain of a nearly asymptotically nonexpansive mapping $T$, then $T$ is asymptotically 
nonexpansive in the intermediate sense. To see this, let $T$ be a nearly asymptotically nonexpansive mapping. 
Then, $$\|T^nx-T^ny\|\le \eta(T^n)(\|x-y\|+a_n)\;\forall\;x,y\in K,\;n\ge 1,$$ which implies that
$$\displaystyle\sup_{x,y\in K}(\|T^nx-T^ny\|-\|x-y\|)\le (\eta(T^n)-1)diam(K)+\eta(T^n)a_n,\;n\ge 1.$$
Hence, $$\displaystyle \limsup_{n\to\infty}\sup_{x,y\in K}(\|T^nx-T^ny\|-\|x-y\|)\le 0.$$
\end{rmk}

\noindent We observe from Remarks \ref{rmk1} and \ref{rmk2}
 that the class of nearly nonexpansive mappings and nearly 
asymptotically nonexpansive mappings are intermidiate classes 
between the class of asymptotically nonexpansive mappings 
and that of asymptotically nonexpansive in the intermediate sense mappings.\\

\noindent Alber \emph{et al.} \cite{Ya} introduced a
more general class of asymptotically nonexpansive mappings called
\emph{total asymptotically nonexpansive mappings} and studied
methods of approximation of fixed points of mappings belonging to
this class.

\begin{definition}\label{def}
A mapping $T:K\to K$ is said to be \emph{total asymptotically
nonexpansive} if there exist nonnegative real sequences
$\{\mu_n\}$ and $\{l_n\},\;n\ge 1$ with ${\mu_n},\;l_n\to 0$ as
$n\to\infty$ and strictly increasing continuous function
$\phi:\mathbb{R^+}\to \mathbb{R^+}$ with $\phi(0)=0$ such that for
all $x, y\in K,$
\begin{eqnarray}\label{e6}
\|T^nx-T^ny\|\le\|x-y\|+\mu_n\phi(\|x-y\|)+l_n,\;n\ge 1.
\end{eqnarray}
\end{definition}

\begin{rmk}
If $\phi(\lambda)=\lambda,$ then \eqref{e6} reduces to
\begin{eqnarray*}
\|T^nx-T^ny\|\le(1+\mu_n)\|x-y\|+l_n,\;n\ge 1.
\end{eqnarray*}
In addition, if $l_n=0$ for all $n\ge 1,$ then total
asymptotically nonexpansive mappings coincide with asymptotically
nonexpansive mappings. If $\mu_n=0$ and $l_n=0$ for all $n\ge 1,$
we obtain from \eqref{e6} the class of mappings that includes the class
of nonexpansive mappings. If $\mu_n=0$ and $l_n= \sigma_n=max\{0,
a_n\},$ where $\displaystyle a_n:=\sup_{x, y\in
K}(\|T^nx-T^ny\|-\|x-y\|)$ for all $n\ge 1,$ then \eqref{e6} reduces to
\eqref{e3} which has been studied as mappings which are asymptotically nonexpansive
in the intermediate sense.\end{rmk}

\begin{rmk} The idea of Definition \ref{def} is to unify various definitions
of classes of mappings associated with the class of asymptotically
nonexpansive mappings and to prove a general convergence theorems
applicable to all these classes of nonlinear mappings.
\end{rmk}

 This work is motivated by the recent paper of Chidume and Ofoedu \cite{Eric}. 
A gap in the proof of Theorem 3.5 in \cite{Eric} is observed. The argument used on page 11, starting from line 8 from bottom to the end of the proof of the theorem is not correct. In this corrigendum, it is our aim to close this gap.

\section {Preliminary}
\noindent In the sequel, we shall need the following
\begin{lemma}\label{lemma1}
Let $\{a_n\},$ $\{\alpha_n\}$ and $\{b_n\}$ be sequences of
nonnegative real numbers such that $$a_{n+1}\le
(1+\alpha_n)a_n+b_n.$$ Suppose that $\displaystyle
\sum_{n=1}^\infty\alpha_n<\infty$ and $\displaystyle
\sum_{n=1}^\infty b_n<\infty.$ Then $\{a_n\}$ is bounded and
$\displaystyle\lim_{n\to\infty}a_n$ exists. Moreover, if in
addition, $\displaystyle\liminf_{n\to\infty}a_n=0,$ then
$\displaystyle\lim_{n\to\infty}a_n=0.$
\end{lemma}

Lemmas:
\begin{lemma}\cite{zeh}\label{lemma310}
Let $E$ be a uniformly convex  Banach space  and $B_R(0)$ be  a closed ball of $E$.  Then there exists  a continuous  strictly increasing convex function $g:[0,\infty)\to [0,\infty)$ with $g(0)=0$
such that
\begin{eqnarray}
   ||\alpha_0 x_0+\alpha_1x_1+\alpha_2x_2+...+\alpha_r x_r||^2 &\leq&  \sum_{i=0}^{r} \alpha_i ||x_i||^2-\alpha_s\alpha_tg(||x_s-x_t||), \nonumber
 \end{eqnarray}
for any  $s,t\in\{0,1,2,...,r\}$  and for $x_i\in B_R(0):=\{x\in E:||x||\leq R\}$, $i=0,1,2,...,r$ with  $\displaystyle \sum_{i=0}^{r}a_i=1$.
\end{lemma}

\section {Main results}
\begin{lemma}\label{uwa}
Let $E$ be a real Banach space, $K$ be a nonempty closed convex
subset of $E$ and $T:K\to K$ be a total
asymptotically nonexpansive mappings with sequences
$\{\mu_{n}\},\;\{l_{n}\}$ $\;n\ge 1$ Suppose that there exist $M,\;M^*>0$ such that
$\phi(\lambda)\le M^*\lambda$ for all $\lambda\ge
M$ then 
\begin{eqnarray*}
\|T^nx-T^ny\|\le (1+\mu_n M^*)\|x-y\|+\mu_n\phi(M)+\ell_n\;\forall\;x,y\in K,\;\forall\;n\ge 1.
\end{eqnarray*}
\end{lemma}
\noindent {\bf Proof.} Since $\phi:[0,+\infty)\to [0,+\infty)$ is strictly increasing continuous function, we have that $\lambda\le M$ implies $\phi(\lambda)\le \phi(M)$; and by the hypothesis, $\phi(\lambda)\le M^*\lambda$ for all $\lambda\ge
M$. It therefore not difficult to see that $\phi(\lambda)\le \phi(M)+M^*\lambda.$ Using this and the fact that $T$ is total asymptotically nonexpansive, we obtain the required inequality. $\Box$\\

\noindent Let $K$ be a nonempty closed convex subset of a real
normed space $E.$ Let $T_1, T_2, ..., T_m:K\to K$ be $m$ total
asymptotically nonexpansive mappings. We first note that the recursion formula (3.1) on page 8 of \cite{Eric} contains a typo. The correct formula is 
\begin{eqnarray}\label{e5}
\displaystyle x_1\in K,\;
x_{n+1}=\alpha_{0n}x_n+\sum_{i=1}^m\alpha_{in}T_i^nx_n,\;n\ge1,
\end{eqnarray}
 where
$\{\alpha_{in}\}_{n\ge 1},\;i=0,1,2,...m$ are sequences in 
$(\gamma_1, \gamma_2)$, for some $\gamma_1,\gamma_2\in (0,1)$ such that $\displaystyle \sum_{i=0}^m\alpha_{in}=1.$\\

\begin{theorem}\label{Th1} (Theorem 3.1 of \cite{Eric})
Let $E$ be a real Banach space, $K$ be a nonempty closed convex
subset of $E$ and $T_i:K\to K,\;i=1, 2, ..., m$ be $m$  total
asymptotically nonexpansive mappings with sequences
$\{\mu_{in}\},\;\{l_{in}\}$ $\;n\ge 1,\;i=1, 2, ..., m$ such that
$\displaystyle F:=\cap_{i=1}^mF(T_i)\ne\emptyset.$ Let $\{x_n\}$
be given by \eqref{e5}. Suppose $\displaystyle\sum_{n=1}^\infty
\mu_{in}<\infty,\;\sum_{n=1}^\infty l_{in}<\infty\;{\rm for\;}i=1, 2, ..., m$
and suppose that there exist $M_i,\;M_i^*>0$ such that
$\phi_i(\lambda_i)\le M_i^*\lambda_i$ for all $\lambda_i\ge
M_i,\;i=1, 2, ..., m,$ then the sequence $\{x_n\}_{n\ge 1}$ is bounded and
$\displaystyle\lim_{n\to\infty}\|x_n-p\|$ exists, $p\in F$.
\end{theorem}

\noindent {\bf Proof.} The proof is exactly the proof of Theorem 3.1 of \cite{Eric}. $\Box$\\

\noindent We now restate and give an alternative proof of Theorem 3.5 of \cite{Eric}.\\

\begin{theorem}(Correctedd version of Theorem 3.5 of \cite{Eric})\label{Th3}
Let $E$ be a uniformly convex real Banach space, $K$ be a nonempty closed convex
subset of $E$ and $T_i:K\to K,\;i=1, 2, ..., m$ be $m$ total asymptotically nonexpansive mappings with
sequences $\{\mu_{in}\},\;\{l_{in}\}\subset[0, \infty)$ such that
$\displaystyle\sum_{n=1}^\infty\mu_{in}<\infty,\;\sum_{n=1}^\infty
l_{in}<\infty,\;i=1, 2, ..., m$ and 
$\displaystyle F:=\underset{i=1}{\overset{m}\bigcap}F(T_i)\ne \emptyset.$ 
 From arbitrary
$x_1\in E,$ define the sequence $\{x_n\}$ by \eqref{e5}. Suppose that
there exist $M_i,\;M_i^*>0$ such that $\phi_i(\lambda_i)\le
M_i^*\lambda_i$ whenever $\lambda_i\ge M_i,\;i=1, 2, ..., m,$ then
$\displaystyle\lim_{n\to\infty}\|x_n-T_i^nx_n\|=0, i=1, 2, ...,
m$. 
\end{theorem}

\noindent {\bf Proof.} Let $p\in F,$ then using the recursion formula \eqref{e5} and Lemma \ref{lemma310}, we have (for any $j\in \{1,2,...,m\}$) that
 \begin{eqnarray}\label{ike}
 \displaystyle \|x_{n+1}-p\|^2&=&\Bigl\|\alpha_{0n}x_n+\sum_{i=1}^m\alpha_{in}T_i^nx_n-p\Bigr\|^2= \Bigl\|\alpha_{0n}(x_n-p)+\sum_{i=1}^m\alpha_{in}(T_i^nx_n-p)\Bigr\|^2\nonumber\\
 &\le &\alpha_{0n}\|(x_n-p)\|^2+\sum_{i=1}^m\alpha_{in}\|T_i^nx_n-p\|^2-\alpha_{0n}\alpha_{jn}g(\|x_n-T_j^nx_n)\nonumber\\ &\le&\alpha_{0n}\|x_n-p\|^2+\sum_{i=1}^m\alpha_{in}\Bigl(\|x_n-p\|+\mu_{in}\phi_i(M_i)\nonumber\\&&+\mu_{in}M_i^*\|x_n-p\|+l_{in}\Bigr)^2-\alpha_{0n}\alpha_{jn}g(\|x_n-T_j^nx_n)\nonumber\\ 
&=&\alpha_{0n}\|x_n-p\|^2+\sum_{i=1}^m\alpha_{in}\|x_n-p\|^2\nonumber\\
&&+\sum_{i=1}^m\alpha_{in}\Bigl(2\|x_n-p\|\big[\mu_{in}\phi_i(M_i)+\mu_{in}M_i^*\|x_n-p\|+l_{in}\big]\nonumber\\&&+\big[\mu_{in}\phi_i(M_i)+\mu_{in}M_i^*\|x_n-p\|+l_{in}\big]^2\Bigr)-\alpha_{0n}\alpha_{jn}g(\|x_n-T_j^nx_n)\nonumber\\ 
&=&\|x_n-p\|^2+\sum_{i=1}^m\alpha_{in}\Bigl(2\|x_n-p\|\big[\mu_{in}\phi_i(M_i)+\mu_{in}M_i^*\|x_n-p\|+l_{in}\big]\nonumber\\&&+\big[\mu_{in}\phi_i(M_i)+\mu_{in}M_i^*\|x_n-p\|+l_{in}\big]^2\Bigr)-\alpha_{0n}\alpha_{jn}g(\|x_n-T_j^nx_n). 
\end{eqnarray} So, since $\alpha_{in}\in (\gamma_1,\gamma_2),\;i=0,1,2,...,m,$ we obtain from \eqref{ike} that
\begin{eqnarray}\label{ike1}\displaystyle \gamma_1^2g(\|x_n-T_j^nx_n)&\le& \alpha_{0n}\alpha_{jn}g(\|x_n-T_j^nx_n)\nonumber\\
&\le&\|x_n-p\|^2-\|x_{n+1}-p\|^2+\sum_{i=1}^m\alpha_{in}\Bigl(2\|x_n-p\|\big[\mu_{in}\phi_i(M_i)\nonumber\\&&+\mu_{in}M_i^*\|x_n-p\|+l_{in}\big]+\big[\mu_{in}\phi_i(M_i)+\mu_{in}M_i^*\|x_n-p\|+l_{in}\big]^2\Bigr)\nonumber\\
&\le&\|x_n-p\|^2-\|x_{n+1}-p\|^2+\gamma_2\sum_{i=1}^m\Bigl(2\|x_n-p\|\big[\mu_{in}\phi_i(M_i)\nonumber\\&&+\mu_{in}M_i^*\|x_n-p\|+l_{in}\big]+\big[\mu_{in}\phi_i(M_i)+\mu_{in}M_i^*\|x_n-p\|+l_{in}\big]^2\Bigr).
\end{eqnarray} But $\mu_{in}\to 0$ and $\ell_{in}\to 0$ as $n\to\infty, i=1,2,...,m$ and by Theorem  \ref{Th1}, $\displaystyle\lim_{n\to\infty}\|x_n-p\|$ exists. Thus, since the summation in \eqref{ike1} is a summation of finite terms, we obtain from \eqref{ike1} that 
 $$\displaystyle \lim_{n\to\infty} g(\|T_j^nx_n-x_n\|)=0\;\forall\;j\in\{1,2,...,m\}.$$ So, properties of the function $g$ (see Lemma \ref{lemma310})  imply that
 \begin{eqnarray}\label{e29}
 \displaystyle \lim_{n\to\infty} \|T_i^nx_n-x_n\|=0,\;i=1,2,...,m. 
 \end{eqnarray} This completes the proof. $\Box$

\begin{theorem}(Corrected version of Theorem 3.6 of \cite{Eric}) 
Let $E$ be a uniformly convex real Banach space, $K$ be a nonempty closed convex
subset of $E$ and $T_i:K\to K,\;i=1, 2, ..., m$ be $m$ continuous total asymptotically nonexpansive mappings with
sequences $\{\mu_{in}\},\;\{l_{in}\}\subset[0, \infty)$ such that
$\displaystyle\sum_{n=1}^\infty\mu_{in}<\infty,\;\sum_{n=1}^\infty
l_{in}<\infty,\;i=1, 2, ..., m$ and 
$\displaystyle F:=\underset{i=1}{\overset{m}\bigcap}F(T_i)\ne \emptyset.$ 
From arbitrary
$x_1\in E,$ define the sequence $\{x_n\}$ by \eqref{e5}. Suppose that
there exist $M_i,\;M_i^*>0$ such that $\phi_i(\lambda_i)\le
M_i^*\lambda_i$ whenever $\lambda_i\ge M_i,\;i=1, 2, ..., m;$ and
that one of $T_1, T_2, ..., T_m$ is compact, then $\{x_n\}$
converges strongly to some $p\in F.$
\end{theorem}

\noindent{\bf Proof.} 
 Observe that from the recursion formula \eqref{e5}, 
\begin{eqnarray}\label{eric2}
\displaystyle \|x_{n+1}-x_n\|&=&\Bigl\|\alpha_{0n}x_n+\sum_{i=1}^m\alpha_{in}T_i^nx_n-x_n\Bigr\|\nonumber\\
&\le &\sum_{i=1}^m\alpha_{in}\|T_i^nx_n-x_n\|\le \gamma_2\sum_{i=1}^m\|T_i^nx_n-x_n\|.
\end{eqnarray}Hence, using \eqref{e29} and  \eqref{eric2}, we obtain that 
\begin{eqnarray}\label{e30}
\displaystyle \lim_{n\to\infty}\|x_{n+1}-x_n\|=0.
\end{eqnarray}
Without loss of generality, let $T_1$ be compact. Since $T_1$ is
continuous and compact, it is completely continuous. Thus, there
exists a subsequence $\{T_1^{n_k}x_{n_k}\}$ of $\{T_1^nx_n\}$ such
that $T_1^{n_k}x_{n_k}\to x^*$ as $k\to\infty$ for some $x^*\in
E.$ Thus $T_1^{n_k+1}x_{n_k}\to T_1x^*$ as $k\to\infty.$ Furthermore, 
\eqref{e29} and the fact that $T_1^{n_k}x_{n_k}\to x^*$ as $k\to\infty$ 
imply that $\displaystyle\lim_{k\to\infty}x_{n_k}=x^*.$
Also from \eqref{e29} $T_2^{n_k}x_{n_k}\to x^*,\;T_3^{n_k}x_{n_k}\to
x^*,\;...,\;T_m^{n_k}x_{n_k}\to x^*$ as $k\to\infty.$ Thus,
$T_2^{n_k+1}x_{n_k}\to T_2x^*,\;T_3^{n_k+1}x_{n_k}\to T_3
x^*,\;...,\;T_m^{n_k+1}x_{n_k}\to T_mx^*$ as $k\to\infty.$ Now,
since from \eqref{e30}, $\|x_{n_k+1}-x_{n_k}\|\to 0$ as $k\to\infty,$ it
follows  that $x_{n_k+1}\to x^*$ as $k\to\infty.$
Next, we show that $x^*\in F.$ Observe that
\begin{eqnarray}\label{e342}
\|x^*-T_ix^*\|&\le&\|x^*-x_{n_k+1}\|+\|x_{n_k+1}-T_i^{n_k+1}x_{n_k+1}\|\nonumber\\
&&+\|T_i^{n_k+1}x_{n_k+1}-T_i^{n_k+1}x_{n_k}\|+\| T_i^{n_k+1}x_{n_k}-T_ix^*\|.
\end{eqnarray}
Taking limit as $k\to\infty$ in \eqref{e342} $\Big($using Lemma \ref{uwa}, \eqref{e29} and \eqref{e30}$\Big)$,
 we have that $x^*=T_ix^*\;(i=1,2,...m)$ and so $x^*\in
F(T_i)\;(i=1,2,...m).$ But by Theorem \ref{Th1},
$\displaystyle\lim_{n\to\infty}\|x_n-p\|$ exists, $p\in F.$ Hence,
$\{x_n\}$ converges strongly to $x^*\in F.$ This completes the
proof. $\Box$\\

\begin{rmk}
We note that the argument (due to the fact that $T$ is uniformly continuous) $\displaystyle \lim_{n\to\infty}\|x_{n+1}-x_n\|=0$ implies $\displaystyle \lim_{n\to 0}\|T^{n+1}x_{n+1}-T^{n+1}x_n\|=0$  as used in \cite{Eric1}, \cite{Eric} and \cite{Udomene} is not always true. To see this, consider $\mathbb{R},$ the set of real numbers endowed with the usual topology and the mapping $T:\mathbb{R}\to \mathbb{R}$ defined by $Tx=3x$ for all $x\in \mathbb{R}.$ It is clear that $T$ is uniformly continuous. Now, let $\{x_n\}_{n\ge 1}$ in $\mathbb{R}$ be a sequence defined by $x_n=1+\frac{1}{n}\;\forall\;n\ge 1.$ We can easily see that $\displaystyle \lim_{n\to\infty}|x_{n+1}-x_n|=0$ but 
\begin{eqnarray*}
|T^{n+1}x_{n+1}-T^{n+1}x_n|&=&\Big|\Big(3^{n+1}+\frac{3^{n+1}}{n+1}\Big)-\Big(3^{n+1}+\frac{3^{n+1}}{n}\Big)\Big|\nonumber\\
&=&\Big|\frac{3^{n+1}}{n+1}-\frac{3^{n+1}}{n}\Big|=\frac{3^{n+1}}{n(n+1)}\to\infty \;{\rm as\;}n\to\infty.
\end{eqnarray*} This follows from the fact that if we define $g:(0,+\infty)\to (0,+\infty)$ by $g(x)=\frac{3^{x+1}}{x(x+1)}$, then by L'Hospital's rule we obtain that 
\begin{eqnarray*}
\displaystyle\lim_{x\to\infty} g(x)=\lim_{x\to\infty}\frac{3^{x+1}}{x^2+x}=3 \ln 3\lim_{x\to\infty}\frac{3^x}{2x+1}=\frac{3(\ln 3)^2}{2}\lim_{x\to\infty} 3^x=+\infty.
\end{eqnarray*}
\end{rmk} Our new method of proof in this corrigendum corrects this error and uniform continuity assumption dispensed with.

\section{Acknowledgement.}
\noindent The autors remain thankful to the reviewer for his/her well articulated comments which helped to improve tha quality of this work.

\end{document}